\documentclass{article}

\usepackage{arxiv}

\usepackage[utf8]{inputenc} % allow utf-8 input
\usepackage[T1]{fontenc}    % use 8-bit T1 fonts
\usepackage{hyperref}       % hyperlinks
\usepackage{url}            % simple URL typesetting
\usepackage{booktabs}       % professional-quality tables
\usepackage{amsthm,amssymb,amsmath,amsfonts}
\usepackage{nicefrac}       % compact symbols for 1/2, etc.
\usepackage{microtype}      % microtypography
\usepackage{lipsum}
\usepackage{graphicx}
\graphicspath{ {./images/} }
\newtheorem{theorem}{\bf Theorem}[section]
\newtheorem{corollary}{\bf Corollary}[section]

\title{Numerical solution of locally loaded Volterra integral equations}

\author{
 Byankin Vladislav \\
  Melentiev Energy Systems Institute\\
  of Siberian Branch of the Russian\\
  Academy of Sciences\\
  Irkutsk,  664033 \\
  \texttt{byankinve@ex.istu.edu} \\
   \And
 Tynda Aleksandr \\
  Penza state University\\
  Penza,  440026 \\
  \texttt{tyndaan@mail.ru} \\
  \And
 Sidorov Denis \\
  Melentiev Energy Systems Institute\\
  of Siberian Branch of the Russian\\
  Academy of Sciences\\
  Irkutsk,  664033 \\
  Harbin Institute of Technology\\
  Harbin, 150001\\
  \texttt{dsidorov@isem.irk.ru} \\
    \And
 Dreglea Aliona\\
Irkutsk National Research Technical University\\
  Irkutsk,  664074 \\
  Harbin Institute of Technology\\
  Harbin, 150001\\
  \texttt{adreglea@ex.istu.edu} \\
}

\begin{document}
\maketitle
\begin{abstract}
Volterra's  integral equations with local and nonlocal loads represent the novel class of integral equations that have attracted considerable attention in recent years. These equations are a generalisation of the classic Volterra integral equations, which were first introduced by Vito Volterra in the late 19th century. The loaded Volterra integral equations are characterised by the presence of a  load which complicates the process of their theoretical and numerical study. Sometimes these equation are called the equations with ``frozen'' argument.
The present work is devoted to the study of Volterra equations with locally loaded integral operators. The existence and uniquness theorems are proved. Among the main contributions is the collocation method for approximate solution of such equations based on the piecewise linear approximation. To confirm the convergence of the method, a number of numerical results for solving model problems are given in the paper.
\end{abstract}

% keywords can be removed

\section{Introduction}

    The term of equations with load was first introduced by  Adam~Nakhushev~\cite{Nakhushev}. The studies of various differential and integral equations  make it possible to  describe more complex processes in science and technologies. Work \cite{CauchyProblem} devoted to the Cauchy problem with a parameter perturbed by a linear functional, which generalises the conventional load. The conditions for the existence of trivial and nontrivial solutions in the neighborhood of the bifurcation point are considered. In the paper \cite{FredholmIntegralEquations} the Fredholm integral equation of the second kind with load is considered. In the paper we constructed methods for regular and irregular cases. Also, the papers have already considered the Hammerstein equation with load in \cite{HammersteinEquations}. In this paper, the equation was studied through a series of nonlinear boundary value problems using the Green's function. In \cite{Kuz-2022}
   the necessary and sufficient conditions for the spectra of the loaded Sturm–Liouville operators, and authors used the term frozen argument instead of the load.

    Let us consider, following the book \cite{Nakhushev}, an integral equation of the form:
    \begin{equation}\label{eqv1:1}
        a_0(t)x(t) + \sum\limits_{j=1}^{m-1}a_j(t)x(t_j) = \lambda \int_{t_0}^t K(t,s) x(s)ds + f(t),\quad t\in \Omega = [t_0,T]
    \end{equation}
    where $t_j, (j=1, 2, \dots ,m-1)$ are defined fixed points of segment $[t_0,T]$, with $t_0<t_2<\dots <t_{m-1}<t_m=T$; $\lambda$  is a constant parameter; $f(t), a_j(t)\in C_{[t_0,T]},j=\overline{0,m-1}$.

    This equation is called the linear one-dimensional loaded Volterra integral equation or the loaded Volterra integral equation (LVIE).

    Let us write the equation (\ref{eqv1:1}) in operator form:
    $$
        Vx(t)=f(t),
    $$
    where the operator $V$ is defined as follows:
    $$
        Vx(t) \equiv a_0(t)x(t) + \sum\limits_{j=1}^{m-1} a_j(t)x(t_j) - \lambda \int_{t_0}^t
        K(t,s)x(s)\,ds.
    $$

    The operator $V$ is a locally loaded operator since it contains the trace of the desired solution in the form of its values $x(t_j), j=1,2,\dots,m-1$, at some points of the set $\Omega$.

  As footnote, let us outline, that to the best of our knowledge, the development of numerical methods for solving equations with loads has not yet been reported in the scientific periodicals.

\section{Existence and uniqueness of solution}

    Let us consider the equation (\ref{eqv1:1}) in the space $C_{[t_0,T]}$, assuming the input functions are continuous:
    \begin{equation}\label{eqv2:1}
        a_0(t)x(t) + \sum\limits_{j=1}^{m-1}a_j(t)x(t_j) + \lambda \int_{t_0}^t K(t,s) x(s)ds = f(t),\quad t\in \Omega = [t_0,T].
    \end{equation}

    Let us also suppose that $a_0(t) \neq 0$. Then, without restricting the generality of the reasoning, let us put $a_0(t) \equiv 1$ in (\ref{eqv2:1}). The kernel $K(s,t)$ in this case has a resolvent $R(t,s,\lambda)$, which can be represented as a uniformly convergent series

    \begin{equation}\label{eqv2:2}
        R(t,s,\lambda) = \lambda K_1(t,s) + \lambda^2K_2(t,s) + \dots + \lambda^nK_n(t,s)
    \end{equation}
    here $K_1(t,s) = K(t,s)$, $K_n(t,s) = \int_s^t K_1(t,z) K_{n-1}(z,s) dz$ and the equation (\ref{eqv2:1}) reduces to the following equation

    \begin{equation}\label{eqv2:3}
        x(t,\lambda) + \sum_{j=1}^m b_j(t,\lambda)x(t_j,\lambda) = F(t,\lambda),
    \end{equation}
    where
    \begin{align*}
        F(t,\lambda) = f(t) + \int_{t_0}^t R(t,s,\lambda) f(s) ds, \\
        b_j (t,\lambda) = a_j(t) + \int_{t_0}^t R(t,s,\lambda) a_j(s) ds.
    \end{align*}
    
    Assuming in (\ref{eqv2:3}) successively $t$ equal to $t_j, j=\overline{1,m}$, to define the load vector $\overline{c} = \big(x(t_1,\lambda), \dots, x(t_m,\lambda) \big)^T$, we obtain a system of linear algebraic equations (SLAE) to determine the load vector $\overline{c} = \big(x(t_1,\lambda), \dots, x(t_m,\lambda) \big)^T$, we obtain a system of linear algebraic equations (SLAE)

    \begin{equation}\label{eqv2:4}
        A(t_j,\lambda) \overline{c} = \overline{d}_{\lambda},
    \end{equation}
    where the matrix $A(t_j,\lambda) \overset{def}{=} [ \delta_{ij} + b_j(t_i,\lambda)]_{ij= \overline{1,m}}$, \(\delta_{ij}\) is--the Kronecker symbol,     \( \overline{d}_{\lambda} = \Big( F(t_j,\lambda), \dots, F(t_m,\lambda) \Big)^T\).

    If $detA \neq 0$, then the SLAE (\ref{eqv2:4}) has a single solution $\overline{C} = A^{-1} \overline{d}_{\lambda}$. In this case, the equation (\ref{eqv2:1}) obviously has a single solution. Hence the theorem follows:

    \begin{theorem}
        Suppose that in equation (\ref{eqv2:1}) all the input functions are continuous and $a_0(t) \equiv 1$. Let $det\,A(t_j,\lambda) \neq 0$. Then the LVIE (\ref{eqv2:1})  has a single solution for all $\lambda$ in the class $C_{[t_0,T]}$. \label{thr1}
    \end{theorem}
    For $\lambda = 0$, the special case $A(t_j,0) = [\delta_{ij} + a_j(t_i)]_{ij= \overline{1,m}}$.

    If $A(t,0)\neq 0$, then by virtue of the continuity of the resolvent on \(\lambda\) and $detA(t,\lambda) \neq 0$ when $|\lambda|detA(t,\lambda) \neq 0$ is small enough, the unambiguous solvability of LVIE  (\ref{eqv2:1}) in $C_{[t_0,T]}$ follows.

    \begin{corollary}
        If $det\,A(t,0) \neq 0$, then if $\lambda$ in equation (\ref{eqv2:1}) is small enough, then the load is uniquely determined and equation (\ref{eqv2:1}) has a unique solution in $C_{[t_0,T]}$ determined from (\ref{eqv2:3}).
    \end{corollary}

    More interesting is the case where the equation contains a value $\lambda$ at which $detA(t,\lambda) = 0$.

    \begin{theorem}
        If $detA(t,0) = 0$ and $rank A(t,\lambda) = r$, then the corresponding homogeneous system will have $m-r$ linearly independent solutions. Correspondingly, the inhomogeneous system (\ref{eqv2:4}) with such $\lambda$ will have no solutions if $\overline{d}(\lambda)$ is not orthogonal to the solution of the conjugate homogeneous system. If $\overline{d}(\lambda)$ is orthogonal to the solution of the conjugate homogeneous SLA, then the inhomogeneous SLA will have a parametric family of solutions depending on $m-r$ free parameters. \label{thr2}
    \end{theorem}

The following theorem follows from this result.

    \begin{theorem}
        Let $det\,A(t,\lambda_0) = 0$, $rank\,A(t,\lambda_0) = r$ and the solvability conditions of the corresponding SLAE (\ref{eqv2:4}) are satisfied, then  integral equation (\ref{eqv2:1}) has a solution dependent on $m-r$ arbitrary constants.
        If the solvability conditions of SLAE (\ref{eqv2:4}) are not fulfilled, then the integral equation (\ref{eqv2:1}) has no solution in $C_{[t_0,T]}$. \label{thr3}
    \end{theorem}

\section{Direct numerical method for the solution of the linear LVIE}

    In this paragraph, we construct a direct numerical method to solve (\ref{eqv1:1}) based on the piecewise linear approximation of the exact solution.

    To discretize equation (\ref{eqv1:1}), we denote by $\Delta_k, k=1,2,\dots,m$, the segments $\Delta=[t_{k-1,t_k}]$, through $\Delta_k, k=1,2,\dots,m$. Then, choosing a small enough sampling step $h$, we cover each of the segments $\Delta_k, k=1,2,\dots,m$, with equally spaced nodes $s_k^l$ defined as follows:

    \begin{equation}\label{eqv3:1}
        s_k^l= t_{k-1} + \frac{(t_k-t_{k-1})}{n_k} \cdot l, \; l=0,1,\dots, n_k-1, \;
        n_k = [\frac{t_k-t_{k-1}}{h}]+1,\; s_m^{n_m}=t_m,
    \end{equation}
    where $[s]$ is the integer part of $s$. The total number of grid nodes (\ref{eqv3:1}) is $N=1+\sum\sum\limits_{k=1}^m n_k$ nodes. For ease of writing, let us re-normalize by $\tau_i, i=0,1,\dots,N,$ all nodes $s_k^l, k=1,2,\dots,m,$ in ascending order.

    Thus, we obtain a mesh of nodes $\tau_i, \; i=0,1,\dots,N,$ consistent with the set of points $t_i, \; j=1,2,\dots,m-1,$ defined in the original equation and satisfying the condition $\underset{j=\overline{1,N}}{\max} (\tau_i - \tau_{j-1}) \leqslant h$.

    Look for an approximate solution $x_N(t)$ of equation (\ref{eqv1:1}) in the form of a piecewise linear function constructed on a grid of nodes (\ref{eqv3:1}):

    \begin{equation}\label{eqv3:2}
        X_N(t)=x(\tau_i-1) + \frac{x(\tau_i) - x(\tau_{i-1})}{\tau_i - \tau_{i-1}} (t-\tau_{i-1}), \; t \in [\tau_{i-1},\tau_i], \; i=\overline{1,N}.
    \end{equation}

    To determine the unknown values $x_i=x(\tau_i),\; i=0,1,\dots,N,$ we require the equation (\ref{eqv1:1}) to be converted to equality at each of the grid points (\ref{eqv3:1}):

    \begin{equation}\label{eqv3:3}
        a_0(\tau_i)x_i + \sum_{j=1}^{m-1} a_j(\tau_i) x_{v_i} = \lambda \int_{t_0}^{\tau_i} K(\tau_i,s) X_N(s)ds+ f(\tau_i), \; i = 0,1,\dots,N,
    \end{equation}
    where $v_i$ is--the number of the grid node $\{ \tau_i \} _{i=0}^N$ coinciding with the value $t_j, \; j = 1,2,\dots, m-1.$ The numbers $v_i$ are determined and fixed during grid generation for a given value of the sampling step $h$.

    Transform (\ref{eqv3:3})

    \begin{equation}\label{eqv3:4}
        a_{0i}x_i + \sum_{j=1}^{m-1} a_{ji} x_{v_i} = \lambda \sum_{p=1}^i \int_{\tau_{p-1}}^{\tau_p} K(\tau_i,s) \Big( x_{p-1} + \frac{x_p-x_{p-1}}{\tau_p-\tau_{p-1}}(s-\tau_{p-1}) \Big) ds+ f_i,
    \end{equation}
    where $a_{ji} = a_j(\tau_i),\; f_i = f(\tau_i),\; j=0,1,\dots,m-1, \;i=0,1,\dots,N$.

    Approximating the integrals in (\ref{eqv3:4}) by the mean rectangle formula, we arrive at a system of linear algebraic equations with respect to the unknown values $x_i,\; i=\overline{0,N}$, of the desired function:

    \begin{equation}\label{eqv3:5}
        a_{0i}x_i + \sum_{j=1}^{m-1} a_{ji} x_{v_i} -f_i = \frac{\lambda}{2} \sum_{p=1}^i (\tau_p-\tau_{p-1}) K \bigg(\tau_i, \frac{\tau_{p-1}+\tau_p}{2} \bigg) (x_{p-1}+x_p),\;i=\overline{0,N}.
    \end{equation}

    For convenience of forming the SLAE matrix by (\ref{eqv3:5}), we introduce the notation:

    \begin{equation*}
        J_p x_{p-1} + J_p x_p=\frac{\lambda}{2} \sum_{p=1}^i (\tau_p - \tau_{p-1}) K\Big(\tau_p, \frac{\tau_{p-1}+\tau_p}{2}\Big)(x_{p-1}+x_p) , \; p = 1, 2, \dots, i.
    \end{equation*}

    Let us also introduce an auxiliary notation for the left part of the SLAE (\ref{eqv3:5}):
    \begin{equation*}
        H_ix_{v_i}=\sum_{j=1}^{m-1} a_{ji} x_{v_i} = (a_{1i} + a_{2i} + \dots + a_{(m-1)i})x_{v_{i}}, \; i=\overline{0,N}.
    \end{equation*}

    Then the system of equations (\ref{eqv3:5}) will take the following form:

    \begin{equation}\label{eqv3:6}
        \begin{cases}
            a_{00}x_0 + H_0 x_{v_0} - f_0 = 0
            \\
            -J_1x_0 + (a_{01}-J_1)x_1 + H_1x_{v_1} - f_1 = 0
            \\
            \cdots\cdots\cdots\cdots\cdots\cdots\cdots\cdots\cdots
            \cdots\cdots\cdots\cdots\cdots\cdots
            \\
            -J_1x_0 + (-J_1-J_2)x_1 + \dots + (-J_{i-2}-J_{i-1})x_{i-1}+ \\
            + (a_{0i} - J_i)x_i + H_ix_{v_i} - f_i = 0.
        \end{cases}
    \end{equation}

    Suppose that the solution of the well-conditioned SLAE (\ref{eqv3:6}) is obtained with sufficient accuracy, for example, by using one of the exact methods, the error of which in this case is limited only by the computational error. Then we can conclude about the following estimation of the error of the proposed numerical method:
    \begin{equation}\label{eqv3:7}
        \varepsilon_h=\|X_N(t)-x(t)\|_{C_{[t_0,T]}} = {\mathcal O} \bigg(\frac{1}{N^2}\bigg) = {\mathcal O}(h^2).
    \end{equation}

{\section{Numerical experiments}

    In order to demonstrate the performance of the proposed numerical method for solving the loaded Volterra integral equations, calculations are performed for two model problems. The algorithm is implemented in C++ programming language, and for solving the system of linear algebraic equations the Gauss--Jordan method with selection of the maximum element is applied.

    \subsection{Model problem 1}
        Consider equation (\ref{eqv1:1}) with the following parameters:
        \begin{gather*}
           t \in [t_0,T] = [0,1];\; m=3; \; t_1=\frac{3}{10}; \;t_2=\frac{1}{2}; \; \lambda = \frac{1}{4}; \\
            a_0(t) = t^2+1; \; a_1(t) = 1-t^3; \; a_2(t) = t-2; \; K(t,s) = t-2s^2; \\
            f(t) = (t^2+1) \cos{t} + (1-t^3) \cos{\bigg(\frac{3}{10}\bigg)} + (t-2) \cos{\bigg( \frac{1}{2} \bigg)}
            + \frac{t^2}{2} \sin{t} - \frac{t}{4} \sin{t} + t \cos{t} - \sin{t}.
        \end{gather*}

        The exact solution in this case is the function \(x(t) = \cos(t)\).

        Table \ref{tab1:1} and Table \ref{tab1:2} show the results of the method for different steps, where the following notations are used: $h$ is the sampling step, \(\varepsilon_h=\|X_N(t)-x(t)\|_{C_{[t_0,T]}}\) is the error,        \( r=\frac{ \ln\frac{\varepsilon_{h_{k-1}}}{\varepsilon_{h_{k}} }}{\ln\frac{h_{k-1}}{h_k}} \) --- order of convergence.
        \newline

        \begin{table}[h]
            \centering
            \begin{tabular}{| c | c | c | c | c | c | c | c | c | c | c | c |}
                \hline
                $h$ & 1/8 & 1/16 & 1/32 & 1/64 & 1/128 & 1/256 \\
                \hline
                $\varepsilon_h$ & 2.20E-04 & 6.99E-05 & 1.91E-05 & 5.04E-06 & 1.30E-06 & 3.30E-07 \\
                \hline
                $r$ & --- & 1.65 & 1.88 & 1.92 & 1.96 & 1.98 \\
                \hline
            \end{tabular}
            \caption{the magnitude of the error and order of convergence.}
            \label{tab1:1}
        \end{table}

        \begin{table}[h]
            \centering
            \begin{tabular}{| c | c | c | c | c | c | c |}
                \hline
                $h$ & 1/512 & 1/1024 & 1/2048 & 1/4096 & 1/8192 & 1/16384\\
                \hline
                $\varepsilon_h$ & 8.31E-08 & 2.09E-08 & 5.23E-09 & 1.31E-09 & 3.27E-10 & 8.18E-11\\
                \hline
                $r$ & 1.99 & 1.99 & 1.93 & 2.07 & 2.00 & 2.00\\
                \hline
            \end{tabular}
            \caption{the magnitude of the error and order of convergence.}
            \label{tab1:2}
        \end{table}

        According to the calculation results, the maximum deviation from the specified function for step $h = 1/8$ is $e = 2.2E-04$. And it is also seen that as the step decreases, the numerical method shows positive convergence.

        Figure \ref{fig1} shows the exact and approximate solution of model problem 1.

        \begin{figure}[h]
            \centering
            \includegraphics[width=0.4\textwidth]{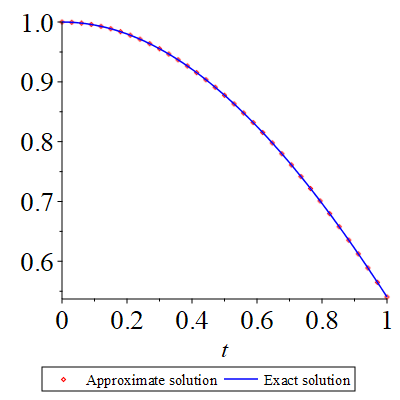}
            \caption{Exact and approximate solution of Model Problem 1, $h = 1/32$}
            \label{fig1}
        \end{figure}

    \subsection{Model problem 2}
        Consider equation (\ref{eqv1:1}) with the following parameters:
        \begin{gather*}
            t \in [t_0,T] = [0,1];\; m=3; \; t_1=\frac{3}{10}; \;t_2=\frac{1}{2}; \; \lambda = \frac{1}{6}; \\
            a_0(t) = \frac{2+t}{3}; \; a_1(t) = t^3-\frac{1}{2}; \; a_2(t) = 2t-t^2; \; K(t,s) = t-2s^2; \\
            f(t) = \Big( \frac{2+t}{3} \Big) e^t + \Big(t^3 - \frac{1}{2} \Big)e^{\frac{3}{10}} + (2t-t^2)e^{\frac{1}{2}} + \frac{e^t \cdot t^2}{3} - \frac{5t \cdot e^t}{6} + \frac{2e^t}{3} + \frac{t}{6} - \frac{2}{3}.
        \end{gather*}

        The exact solution in this case is the function \( x(t) = e^t \)

        Table \ref{tab2:1} and Table \ref{tab2:2} show the results of the method for different steps, where the following notations are used: $h$ is the sampling step, \(\varepsilon_h=\|X_N(t)-x(t)\|_{C_{[t_0,T]}}\) is the error,        \( r=\frac{ \ln\frac{\varepsilon_{h_{k-1}}}{\varepsilon_{h_{k}} }}{\ln\frac{h_{k-1}}{h_k}} \) is order of convergence.

        \begin{table}[h]
            \centering
            \begin{tabular}{| c | c | c | c | c | c | c |}
            \hline
                h & 1/8 & 1/16 & 1/32 & 1/64 & 1/128 & 1/256\\
                \hline
                e & 7.99E-04 & 2.33E-04 & 6.80E-05 & 1.85E-05 & 4.73E-06 & 1.19E-06\\
                r & & 1.78 & 1.78 & 1.88 & 1.97 & 1.98\\
                \hline
            \end{tabular}
            \caption{the magnitude of the error and order of convergence.}
            \label{tab2:1}
        \end{table}

        \begin{table}[h]
            \centering
            \begin{tabular}{| c | c | c | c | c | c | c |}
            \hline
            h & 1/512 & 1/1024 & 1/2048 & 1/4096 & 1/8192 & 1/16384\\
            \hline
            e & 3.02E-07 & 7.59E-08 & 1.90E-08 & 4.75E-09 & 1.19E-9 & 2.97E-10\\
            r & 1.98 & 1.99 & 1.93 & 2.07 & 2.00 & 2.00\\
            \hline
        \end{tabular}
            \caption{the magnitude of the error and order of convergence.}
            \label{tab2:2}
        \end{table}

        According to the calculation results, the maximum deviation from the specified function for step $h = 1/8$ is $e = 7.99E-04$. It is also seen that as the step decreases, the numerical method shows positive convergence.

        Figure \ref{fig2} shows the exact and approximate solution of model problem 2.

        \begin{figure}[h]
            \centering
            \includegraphics[width=0.4\textwidth]{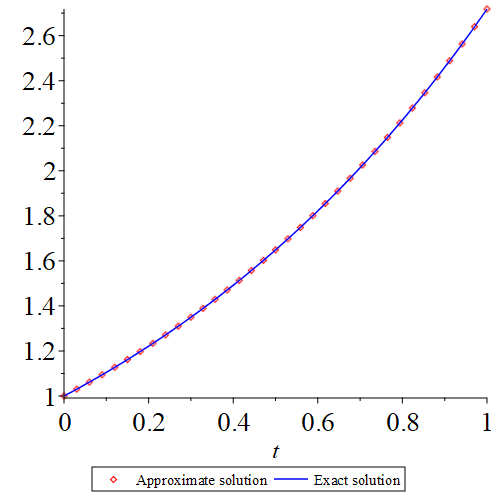}
            \caption{Exact and approximate solution of Model Problem 2, $h = 1/32$}
            \label{fig2}
        \end{figure}

\section{Conclusion}

    In this paper, the application of the new numerical method of piecewise linear approximation for the linear loaded Volterra integral equations of the second kind is considered. Efficiency of the novel numerical method was carried out using two model problems. The numerical solution of the nonlinear Volterra and Hammerstein integral equations will be discussed in the following papers.

\section{Acknowledgement}
This work was supported by the Ministry Science and Higher Education of Russian Federation, project No. FZZS-2024-0003.

\bibliographystyle{unsrt}


\begin{thebibliography}{}

\end{thebibliography}


\begin{thebibliography}{1}

\bibitem{Nakhushev}
Nakhushev A.M.
\newblock Loaded equations and their application.
\newblock {\em Research Institute of Invitational Mathematics and Automation KBNTs RAS}, Moscow: Nauka, 2012. 232 p.

\bibitem{CauchyProblem}
 Dreglea Sidorov L. R.,  Sidorov N. A.
\newblock Cauchy problem with a parameter perturbed by a linear functional.
\newblock {\em Proceedings of Science and Technology, Modern Math. and its Appl. Thematic review.}, 2024, vol. 237, pp. 10–17.

\bibitem{FredholmIntegralEquations} Dreglea Sidorov L. R.,  Sidorov N.,  Sidorov D.
\newblock The linear Fredholm integral equations with functionals and parameters.
\newblock {\em Bul. Acad. Ştiinţe Repub. Mold. Mat.}, 2023, vol. 2, pp. 83–91.

\bibitem{HammersteinEquations}Sidorov N.A., Dreglea Sidorov L.R.D.
\newblock On the Solution of Hammerstein Integral Equations with Loads and Bifurcation Parameters.
\newblock {\em The Bulletin of Irkutsk State University. Series Mathematics}, 2023, vol. 43, pp. 78–90.

\bibitem{ISU-2022}	Tynda A., Noeiaghdam S., Sidorov D.
\newblock Polynomial spline collocation method for solving weakly regular Volterra integral equations of the first kind.
\newblock {\em The Bulletin of Irkutsk State University. Series: Mathematics.}, 2022. vol. 39. pp. 62-79.

\bibitem{Kuz-2022}	
Kuznetsova M.
\newblock Necessary and sufficient conditions for the spectra of the Sturm–Liouville operators with frozen argument.
\newblock{\em Applied Mathematics Letters.},  2022. vol. 131. pp. 108035.

\end{thebibliography}

\end{document}